\documentclass[11pt]{article}
\usepackage{amsmath}
\usepackage{amssymb}
\begin{document}
\title{An Operator Theoretical  Proof Of The Fundamental Theorem  Of Algebra }
\author{Ali Taghavi \\Department of  Mathematics, Damghan University of  Basic Science,Iran\\
alitghv@yahoo.com,taghavi@dubs.ac.ir}

\maketitle
\begin{abstract}
We  give a proof of the fundamental theorem of algebra with operator
theoretical approach
\end{abstract}

\textbf{Introduction}:
\\

 The fundamental theorem of algebra states every polynomial
 equation $P(z)=0$ with complex coefficients has at least one
 complex root.Various methods have been used for a proof of this
 historical theorem . Topological methods,complex variable methods(using Liouville theorem) and algebraic methods are examples of such arguments.
 In spite of of its algebraic nature,almost all proofs of this
 theorem
 involve topology or analysis.
 In this note,we  present a proof with an operator theoretical view point,

 First we present some preliminaries :
 \\

 \textbf{Preliminaries}:
 \\
 A Banach space is a vector space $X$ over the real or
complex numbers with a norm $|| . ||$ such that every Cauchy
sequence
 (with respect to the metric $d(x, y) = ||x - y||$) in $X$ has a limit
 in $X$.A Hilbert space,is a vector space $H$ with an inner product
 $< . >$,such that H with the norm which induced by $< .>$ is a
 Banach space.
 \\

 Assume  $X$ is a Banach space,a continuous linear transformation from $X$ to $X$ is called a bounded
 operator.Here the word $bounded$  signifies that any such   operator
 maps
 the unit disk of $X$ to a bounded set of $X$.
 We denote by $B(X)$, the linear space of all
 bounded operators on $X$.Consider the standard norm on $B(X)$ as
 follows:for a  given S  in  $B(X)$,define
 \\

  $||S||=\sup_{|z|=1}|S(z)|$.This norm induce a topology on
 $B(X)$,which is called,the norm topology.
 \\
 By an isometry on $X$ we mean a linear operator on $X$ which
 preserves the norm of $X$.
 \\

 $l^{2}$ is the Hilbert space  of all sequences  of complex numbers
 $(a_{i})$ with convergence series $\sum|a_{i}|^2$. $l^{2}$ is equipped with
 the inner product $\sum_{n=1}^{\infty}a_{n}\bar{b_{n}}$
  which induce
 the norm,$|(a_{n})|_{2}=\sqrt{\sum_{i=1}^{\infty}|a_{i}|^{2}}$.We say
 $(a_{n})$ is an eventually zero sequence if there exist $k$ such
 that $a_{i}=0$ for all $i>k$
 \\

 By $i-shift$ operator on $l^{2}$ ,we
mean the operator which sends $(a_{i})$ to
(0,0,\ldots$,a_{1},a_{2},\ldots$),namely this operator adds $i$
zeroes at the first terms of the sequence $(a_{n})$.
 \\
 \\

 $Hol(\mathbb{C})$ is the space of all entire maps from $\mathbb{C}$
 to $\mathbb{C}$.\\
 Using Taylor expansion, we embed $Hol(\mathbb{C})$ in $l^2$ as a dense linear subspace, as follows:

 Let f be an entire map with Taylor expansion $f(z)=\sum
 a_{n}z^{n}$,then the sequens $(a_{n})$ is an element of $l^2$,because
 convergence   of $\sum a_{n}z^{n}$ for all $z$ with  $|z|>1$ implies that the
 series is absolutely converge for $z=1$,then $\sum|a_{n}|$
 converges,and this implies that $\sum|a_{n}|^2$ converges, too.
\\
 Note that
 this embedding sends  a polynomial $P(z)=a_{n}z^{n}+a_{n-1}z^{n-1}+\ldots+a_{0}$
 to the sequence $(a_{0},a_{1},\ldots,a_{n},0,0,\ldots)$,this
 sequens  is eventually zero.As we prove in the following lemma 1,the space of all sequence which are
 eventually zero is a dense linear subspace of $l^2$.So we consider
 $Hol(\mathbb{C})$ as a dense linear subspace of $l^2$,since
 $Hol(\mathbb{C})$ contains all polynomials in one variable.
 \\

  \textbf{Lemma 1}: The space of all eventually zero sequences is a
  dense linear subspace of  $l^{2}$.
  \\

  \textbf{Proof of lemma 1}:Assume  $(a_{n})$ is an element of
  $l^{2}$,and $\varepsilon$ is given,we give an eventually zero
  sequence $(b_{n})$ such that $|(a_{n})-(b_{n})|_{2}\leq
  \varepsilon$.With  the Cauchy criteria for  the series $\sum|a_{n}|^{2}$
  ,we have an integer k such that
  $\sum_{i=k+1}^{\infty}|a_{i}|^{2}<\varepsilon^{2}$, now put
  $(b_{n})$=$(a_{1},a_{2},\ldots,a_{k},0,0,\ldots)$,then\\
  $|(a_{n})-(b_{n})|_{2}=\sqrt{\sum_{i=k+1}^{\infty}|a_{i}|^{2}}\leq
  \varepsilon$.
\\
\\

 We assign to a polynomial $P(z)$, a bounded operator on
$l^2$
 which restriction to $Hol(\mathbb{C})$ is equal to the map of  multiplication
 by $P(z)$,i.e.the map which sends $f$ to $Pf$.\\To prove this,it suffices  to assign a
 bounded operator on $l^{2}$ correspond to monomial  $z^{i}$,since any polynomial is a
 linear combination of  monomials $z^{i}$.
 \\
 \\
 Assume $(a_{n})$ is the coefficients  of Taylor series of an entire map f(z),then the
 coefficients of Taylor expansion of $z^{i}f(z)$ is
 $(0,0\ldots,0,a_{1},a_{2},\ldots)$,where $a_{1}$ is placed at $i+1$-th term .So multiplication by $z^{i}$ as a linear operator on $Hol(\mathbb{C})$
 can be extended to $i-shift$ operator on  $l^{2}$.
\\

Let $H$ be a Hilbert space,and $T$ be a bounded operator on
$H$,which
 range is a closed subspace of $H$,we say T is a Fredholm  operator
 if both kernel and ko-kernel of $T$ are finite dimensional linear
 space where co-kernel of $T$ is the quotient   space $H/rang(T)$\\
By definition,the Fredholm index of $T$ is dim of kernel of $T$
minus dim of  co-kernel of $T$.
\\
In the following we state,without proofs, some essential properties
of Fredholm operators,(for more information about Fredholm operator
theory see $[1]$)or $[2]$:
\\

 The space of  all Fredholm operators on $H$ is denoted by $Fred(H)$.
 $Fred(H)$ is an open subset of the space of all bounded operators and index is a
 continuous map from $Fred(H)$ to integers, $\mathbb{Z}$.Thus if $T$  is a Fredholm operator,there is a neighborhood of
 T,in $B(H)$,with the norm topology,
 which elements are Fredholm operators of the same index as index of $T$.In fact Fredholm index is a continuous map which is locally
 constant.This property is called,"Invariance of Fredholm index with small perturbation".
 \\

In the next part, we present a complete proof for the Fundamental
theorem of
 algebra.
 \\

\textbf{The proof }\\
Assume $P(z)=z^{n}+a_{n-1}z^{n-1}+\ldots$ has no root,then
$Q(z)=\varepsilon^nP(z/\varepsilon^n)$ does not have any root
,too.we have $Q(z)=z^{n}+\epsilon
a_{n-1}z^{n-1}+\epsilon^{2}a_{n-2}z^{n-2}+\ldots+\epsilon^{n}a_{0}$

 We define a  bijective operator on $Hol(\mathbb{C})$ by multiplication by $Q(z)$,since
$Q(z)$has no root.This operator can be extended to a bounded linear
operator on $l^{2}$,as we explained above.We call this
extension,Q,again\\

  Now $Q(z)$ is the perturbation of $n-shift$ operator,that is,$Q$ is sufficiently close to
 $n-shift$ operator if $\varepsilon$ is sufficiently small.
  Further it can be easily proved that  $n-shift$ operator satisfies in the following two conditions :

 1)The  Fredholm  index of $n-shift$ operator is $-n$\\

 2)$n-shift$ operator is an isometry
 \\

 Since $Q$ is a perturbation of $n-shift$ operator ,it is a Fredholm operator of index $-n$,this is because of invariance of fredholm index
 with small perturbation.Further we use the following lemma 2 to prove that Q is a one to one,Fredholm operator of index $-n$ and satisfies the inequality
$|Q(z)|>k|z|$ for some constant $k$.

  Obviously any  $Q$ with such properties can not be a surjective operator because the kernel of $Q$ has zero dimension,so the codimension of
  the range of $Q$ is $n$.The fact that $Q$ is not surjective contradicts to the
  following lemma 3,and this would complete the proof of the fundamental theorem of algebra.
  \\

\textbf{Lemma 2} :Let  H be  a Banach space and T be an isometry on
H,then there is a neighborhood W of T,in the space of bounded
operator on H,such that for every S in W ,we have $|S(z)>k|z|$,for a
constant k depend on S.
\\

\textbf{Proof of lemma 2}:Assume $|T-S|\leq1/2$ ,then the inequality
$|z|=|T(z)|\leq |T(z)-S(z)|+|S(z)|$ implies $|S(z)|\geq |z|/2$.
\\

\textbf{Lemma 3} :Let  $F$ be a dense subspace of a Banach space
 $E$,and T is a bounded linear operator on $E$,which maps F onto F,further $|Tz)|>k|z|$ for some $k$,then
 $T$ is a surjective operator onto E.
 \\

 \textbf{Proof of lemma 3}: If $T$ is not surjective ,then there is an element $e$  of $E$ which is
 not in the range of $T$.Since $F$ is dense ,there is a sequence $f_{n}$ of elements
 of $F$ which converge to $e$ and $f_{n}=T(b_{n})$ for some $b_{n}$ in $F$.The condition $|Tz)|>k|z|$
 implies that the pre-image of a Cauchy sequence under $T$ is a Cauchy sequence .So $b_{n}$ is a Cauchy sequence,
 then converges to some $b_{*}$ in $E$.
 From  continuity of $T$ we have $T(b_{*})=e$,that is $e$ is in the
 image of $T$,which contradicts to our first assumption.This
 completes the proof of lemma 3
 \\


\begin{thebibliography}{99}

\bibitem{1}

\emph{John B. Conway,A Course In Functional
Analysis,Springer-Verlag,1985},

\bibitem{2}

\emph{John B. Conway,A Course In Operator Theory,American
Mathematical Society, Providence, RI, 2000}


\end{thebibliography}
\end{document}